\documentclass{article} \usepackage{amssymb,amsmath,eucal}

\begin{document}

\def\R {{\mathbb R }}
 \def\C {{\mathbb C }}
  \def\Z{{\mathbb Z}} 
  \def\H{{\mathbb H}}

\def\vr{\mathbf {wr}}

\newcommand{\arcsh}{\mathop{\rm arcsh}\nolimits}
 \newcommand{\sh}{\mathop{\rm sh}\nolimits} 
\newcommand{\ch}{\mathop{\rm ch}\nolimits}

\def\SL{{\rm SL}}

\def\U{{\rm U}}
\def\O{{\rm O}} 
\def\Sp{{\rm Sp}} 
\def\SO{{\rm SO}}

\def\ov{\overline} 
\def\phi{\varphi} 
\def\epsilon{\varepsilon}
\def\kappa{\varkappa}

\def\le{\leqslant} 
\def\ge{\geqslant}

\def\F{\,\,{}_2F_1} 
\def\FF{\,\,{}_3F_2}

\renewcommand{\Re}{\mathop{\rm Re}\nolimits} 
\renewcommand{\Im}{\mathop{\rm Im}\nolimits}

\newcounter{sec} \renewcommand{\theequation}{\arabic{sec}.\arabic{equation}}

\newcounter{fact} \def\fact{\addtocounter{fact}{1}{\scc \arabic{fact}}}

\newcounter{punkt}
\def\punkt{\addtocounter{punkt}{1}{\arabic{sec}.\arabic{punkt}.  }}

\begin{center} \Large

\bf Some continuous analogs of expansion
\linebreak
in Jacobi polynomials and
\linebreak
vector valued hypergeometric orthogonal bases

\medskip

\large
\sc Neretin Yu.A.

\end{center}

{\sc Abstract.} {\small
We write spectral decomposition of the hypergeometric differential operator on the contour $\Re z=1/2$
(multiplicity of spectrum is 2). As a result, we obtain an integral transform that differs from the Jacobi (or Olevsky) transform. We also try to do a step towards vector-valued special functions and construct a
${}_3F_2$-orthogonal basis in the space of functions having values in 2-dimensional space. This basis is lying in an analytic continuation of continuous dual Hahn polynomials with respect to number $n$ of a polynomial.

 In Addendum, we discuss  $\C^2$-value analogs
 of the Meixner--Pollachek orthogonal system and
 also perturbations of Laguerre, Meixner, and
 Jacobi polynomials.
 }

\addtocounter{sec}{1}

\bigskip

 {\bf\large $\mathbb \S$1.  Representation-theoretical motivation
 and formulation of
results.}

\bigskip

{\bf \punkt Continuous analogs of expansion in Jacobi polynomials.}  The
present work is an counterpart on the level
of special functions
 of Molchanov's paper \cite{Mol-tenzor} on tensor
products of unitary representations of the group $\SL_2(\R)$.

In the classical analysis, there are well-known  the expansion in Jacobi polynomials
and also its continuous analogue.  The latter is the index hypergeometric
transform of H.Weyl, \cite{Wey} (it is also called by Olevsky transform
\cite{Ole},
Jacobi transform, generalized Fourier transform, see Koornwinder's
survey \cite{Koo2}, see also \cite{FK}, \cite{Koo1}, \cite{Ner-index}).  These
classical constructions have a transparent and important representation
theoretic interpretation (see  \cite{VK}).  Spherical functions of the
projective spaces
\begin{equation}
\!\!\!\!\!\!  \O(n+1)/\O(n)\times \O(1),
\,\,
\quad \U(n+1)/\U(n)\times \U(1),\,\,%
\quad \Sp(n+1)/\Sp(n)\times \Sp(1)
\end{equation}
over $\R$, $\C$ and the quaternion field $\H$ are  the Jacobi
polynomials $P_m^{\alpha,\beta}$ for some special values of the parameters
$\alpha$, $\beta$.  The theorem about
 decomposition of $L^2$ on these spaces into
a direct sum of irreducible representations is a corollary of the theorem on
expansion of a function into a series in Jacobi polynomials.

For hyperbolic spaces
\begin{equation} \O(n,1)/\O(n)\times \O(1),\,\,%
 \qquad
\U(n,1)/\U(n)\times \U(1),\,\,%
 \qquad \Sp(n,1)/\Sp(n)\times \Sp(1)
\end{equation}
the analogy of an expansion in the Jacobi polynomials is the {\it index
hypergeometric transform} (see  \cite{Ole}, \cite{Koo2})
\begin{equation}
g(s)= \frac 1{\Gamma(b+c)}\int_0^\infty f(x) \F(b+is,b-is;b+c;-x)
x^{b+c-1} (1+x)^{b-c} dx
.\end{equation}
The problem of decomposition of $L^2$
on these spaces is reduced to the inversion formula
$$
f(x)= \frac
1{\pi\Gamma(b+c)}\int_0^\infty g(s) \F(b+is,b-is;b+c;-x) \Bigl|\frac{\Gamma(b+is)
\Gamma(c+is)}{\Gamma(2is)}\Bigr|^2 ds
$$
 for this integral transform.

Nevertheless, these two classical constructions
(i.e., the expansion in the Jacobi
polynomials and the index hypergeometric transform
 are not sufficient in the analysis on
pseudo-Riemannian symmetric spaces of rank 1 (this class of spaces includes, in
particular, other real forms of the spaces (1.1)--(1.2)).  This force to think
that there exists some another analog (or analogs) of the index hypergeometric
transform (1.3).  We construct one such transform%
\footnote{Molchanov
\cite{Mol-tenzor}, \cite{Mol-rank} \cite{Mol-VINITI}
 uses another method.  Below
we introduce the hypergeometric
 differential operator (2.1).
  It arises if we restrict the Laplace
operator in tensor product of
two unitary representations of
the universal covering group $\SL_2(\R)^\sim$
 of $\SL_2(\R)$ to
eigenspaces of the rotation group $\SO(2)$.
Molchanov  \cite{Mol-tenzor} (he consider
only the group ${\rm PSL}_2(\R)$) restricts the
Laplace operator to functions that are invariant with respect to the group of
diagonal matrices.  As a result, he obtains
 the Legendre differential operator on
some contour containing  singular points of the Legendre equation.}
\footnote{It is interesting to understand, is it possible to give another proof
of the Molchanov's Plancherel formula \cite{Mol-rank}
 for rank 1 symmetric spaces
$G/H$ using the classical index hypergeometric transform (in the extended variant
\cite{DS}, XIII or \cite{FK}) and the double index hypergeometric transform
constructed below in 1.4.  Precisely, let $K$ be a maximal compact subgroup in
$G$ , let $V$ be an irreducible $K$-module.  Consider the action of the Laplace
operator on $\mathrm{Hom}_K(V,L^2(G/H))$.  Is it correct, that in all the cases
we will obtain the hypergeometric operator (2.3) on the segment $[0,1]$ (the
classical case) or on the contour $\Re z=1/2$ (our case)?}.

\smallskip

{\bf\punkt The problem on vector-valued bases.}  The Askey--Wilson hierarchy of
hypergeometric orthogonal polynomials is well known, see for instance, \cite{KS}, or \cite{AAR},
Chapter 6.  Almost all these polynomials (probably all) appear in the
representation theory of the group $\SL_2(\R)$ (see, for instance, \cite{VK},
\cite{Ros1}--\cite{Ros2}).

Consider the tensor product of a unitary highest weight representation of
$\SL_2(\R)$ and a lowest weight representation.  In each factor, there is a
canonical orthogonal basis consisting of $\SO(2)$-eigenfunctions.  Hence there
exists a canonical basis in the tensor product.  This basis consists of
continuous dual Hahn  polynomials
(\cite{Zha}).  Recall (see \cite{AAR}, 6.10, \cite{KS}) that
they are the polynomials  $p_n(s^2)$ orthogonal with respect to the measure
$$
\Bigl|\frac{
\Gamma(a+is)\Gamma(b+is) \Gamma(c+is)}{\Gamma(2is)}
\Bigr|^2 ds
$$
 on the line
$s\in\R$;  the explicit formula is
\begin{equation}
p_n(s^2)=(a+b)_n(a+c)_n \FF\Bigl[\begin{matrix} -n,a+is,a-is\\a+b,a+c
\end{matrix};1\Bigr]
.  \end{equation}
 But the same problem has sense for each
pair of unitary representations of $\SL_2(\R)$
or its universal covering.  These tensor products
(generally) have multiplicity 2 (see \cite{Puk}, \cite{Mol-tenzor}).  Hence, we
must obtain some orthogonal bases consisting of $\C^2$-valued functions.

\smallskip

{\bf \punkt Notations.}  We use the standard notations for the Pochhammer symbol
$$(a)_n=a(a+1)\dots(a+n-1)$$
 and for hypergeometric functions
 \begin{align*}
\F[a,b;c;z]:&= \sum_{n=0}^\infty \frac{(a)_n(b)_n}{(c)_n n!  } z^n;\\ \FF
\Bigl[\begin{matrix} a,b,c\\d,e \end{matrix};z\Bigr] :&= \sum_{n=0}^\infty
\frac{(a)_n(b)_n (c)_n}{(d)_n (e)_n n!  } z^n 
.\end{align*}

{\bf \punkt Double index hypergeometric transform.}  Fix $0\le \alpha\le 1/2$,
$\beta\in \R$.  Assume that $\alpha+i\beta\ne 0$.

Consider the space of $\C^2$-valued functions of the half-line $s>0$.  An
element of this space can be considered as a pair of scalar-valued functions
$(\phi_1(s),\phi_2(s))$.  Let us introduce the scalar product in this space by
the formula 
\begin{multline*}
\bigl\langle(\phi_1,\phi_2),(\psi_1,\psi_2)\bigr\rangle = \frac1{2\pi} 
\int\limits_0^\infty \Bigl[ r_{11}(s) \phi_1(s) \ov{\psi_1(s)} + r_{12}(s)
\phi_1(s) \ov{\psi_2(s)} + \\+ r_{21}(s) \phi_2(s) \ov{\psi_1(s)} + r_{22}(s)
\phi_2(s) \ov{\psi_2(s)} \Bigr] \frac{ds}{|\Gamma(2is)|^2} 
,\end{multline*}
where $r_{ij}(s)$ are given
by
 \begin{multline}
 R(s)=\begin{pmatrix}
r_{11}(s)&r_{12}(s)\\ r_{21}(s)&r_{22}(s) \end{pmatrix} :=\\:= \begin{pmatrix}
\Gamma(\tfrac 12-\alpha-is) \Gamma(\tfrac 12-\alpha+is) & \Gamma(\tfrac
12-i\beta-is) \Gamma(\tfrac 12-i\beta+is) \\ \Gamma(\tfrac 12+i\beta-is)
\Gamma(\tfrac 12+i\beta+is) & \Gamma(\tfrac 12+\alpha-is) \Gamma(\tfrac
12+\alpha+is) \end{pmatrix} 
.\end{multline} 
It is convenient to write this
scalar product in the vector form 
$$
\bigl\langle(\phi_1,\phi_2),(\psi_1,\psi_2)\bigr\rangle =
\frac1{2\pi}\int_0^\infty \begin{pmatrix}\phi_1(s) &\phi_2(s) \end{pmatrix} R
(s) \begin{pmatrix}\ov{\psi_1(s)} \\ \ov{\psi_2(s)} \end{pmatrix}
\frac{ds}{|\Gamma(2is)|^2}
.$$ 
Thus we obtain the Hilbert space of $\C^2$-valued
functions.  We denote it by $H_{\alpha,\beta}$.

Let $x\in \R$, $s\ge 0$.  Consider two functions $Q_1(\alpha,\beta; x,s)$,
$Q_2(\alpha,\beta; x,s)$ given by the formulae%
\footnote{$Q_1$, $Q_2$
are almost $L^2$-eigenfunctions of the differential
operator $D$ defined below (2.1).}
\begin{multline}
Q_1(\alpha,\beta; x,s) =\frac 1{\Gamma(\alpha+ i\beta)}
(\tfrac12+ix)^{-(\alpha+i\beta)/2} (\tfrac12-ix)^{(\alpha+i\beta)/2-is-1/2}
\times \\ \times \F \Bigl[\begin{matrix} 1/2-\alpha+is,1/2-i\beta+is\\
1-\alpha-i\beta \end{matrix}\,;\frac{ix+1/2}{ix-1/2} \Bigr] 
;\end{multline}
 $$
\!\!\!\!\!\!\!\!\!\!\!\!\!\!\!\!\!\!\!\!\!\!\!\!\!\!\!\!
\!\!\!\!\!\!\!\!\!\!\!\!\!\!\!\!\!\!\!\!\!\!\!\!\!\!\!\!
\!\!\!\!\!\!\!\!\!\!\!\!\!\!\!\!\!\!\!\!\!\!\!\!\!\!\!\!
\!\!\!\!\!\!\!\!\!\!\!\!\!\!\!\!\!\!\!\!\!\!\!\!\!\!\!\!  Q_2(\alpha,\beta; x,s)
= Q_1(-\alpha,-\beta; x,s) 
.$$

For a function $f\in L^2(\R)$, we define the pair of functions $\phi_1(s)$,
$\phi_2(s)$ by
 \begin{equation}
  \phi_j(s):=\int_{-\infty}^\infty
f(x)\ov{Q_j(\alpha,\beta;x,s)}\, dx, \qquad j=1,2
 .\end{equation}

{\sc Theorem 1.1.}  {\it Let $0\le\alpha\le 1/2$, $\alpha+i\beta\ne0$.}

a) {\it The operator $f\mapsto(\phi_1,\phi_2)$ is a unitary bijective operator from
$L^2(\R)$ to $H_{\alpha,\beta}$.}

\smallskip

b) {\it The inversion formula is}
 \begin{equation} 
 f(x):=\frac
1{2\pi}\int_0^\infty \begin{pmatrix} Q_1(\alpha,\beta; x,s) & Q_2(\alpha,\beta;
x,s) \end{pmatrix} R(s) \begin{pmatrix} \phi_1(s)\\ \phi_2(s)\end{pmatrix}
\frac{ds}{|\Gamma(2is)|^2}
 .\end{equation}

{\sc Remark.} Let $A(s)$ be a
$2\times2$-matrix valued
function on $\R$. Any
 calibration
$$
(Q^\circ_1(s),Q^\circ_2(s)):=
(Q_1(s),Q_2(s))A(s)^{-1},\qquad
R^\circ(s):=A(s)^*R(s)A(s)
$$
 gives another
form of our integral transformation.
It is interesting to find matrices
$A(s)$, for which
the matrix $R^\circ(s)$ and the vector
$(Q^\circ_1(s),Q^\circ_2(s))$
remain relatively simple.

\smallskip

{\bf \punkt Generalization for $\alpha>1/2$.}  Now, let $\alpha> 1/2$,
$\beta\in\R$, and $\alpha+ i\beta\ne 0,1,2,\dots$.  Denote by $n$ the integral
part of $\alpha-1/2$.  Consider the finite dimensional linear space
$W_{\alpha,\beta}$, consisting of vectors $(c_0,c_1,\dots,c_n)$.  The scalar
product in $W_{\alpha,\beta}$ is given by
$$
 \langle c,c'\rangle=\frac
1{2\pi}
\sum_{k=0}^n\frac{2\alpha-2k-1}{\Gamma(2\alpha-k)k!}  c_k {\ov c_k}'
.$$

Next, define the functions%
\footnote{These functions are $L^2$-eigenfunctions
of the differential operator $D$ defined below}
\begin{multline}
R(\alpha,\beta;x;k):=\\
:=\frac{1}{\Gamma(\alpha+i\beta)}
(\tfrac12+ix)^{-(\alpha+i\beta)/2} (\tfrac12-ix)^{-(\alpha-i\beta)/2}
\F\Bigl[\begin{matrix}-k,k-2\alpha+1\\1-\alpha-i\beta\end{matrix};
\tfrac12+ix\Bigr] 
.\end{multline}

Now, consider the linear operator 
\begin{equation} 
J_{\alpha,\beta}:L^2(\R)\to
H_{\alpha,\beta} \oplus W_{\alpha,\beta}
 ,\end{equation} 
 given by
$$
f\mapsto(\phi_1,\phi_2, \theta)
,$$
where $\phi_1$, $\phi_2$ is defined by
(1.7) as above, and the coordinates of the vector $\theta\in W_{\alpha,\beta}$
have the form 
$$
\theta_k=\int_{-\infty}^\infty f(x)\ov{R(\alpha,\beta;x;k)}\,dx
.$$

{\sc Theorem 1.2.} {\it The operator $J_{\alpha,\beta}$
is a unitary bijective operator.}

\smallskip

{\bf \punkt Romanovski polynomials.}  The Romanovski polynomials \cite{Rom} are
the polynomials on $\R$ orthogonal with respect to the weight
 $$
(\tfrac12+ix)^{-(\alpha+i\beta)} (\tfrac12-ix)^{-(\alpha-i\beta)} \,dx
.$$
This
weight decreases as $x^{-2\alpha}$, and hence it is possible to orthogonalize
only finite number of power functions $1$, $x$, $x^2$, \dots.  Romanovski
polynomials are defined by the formula
 $$
 \F[-k,k-2\alpha+1;1-\alpha-i\beta;
\tfrac12+ix] 
,$$ 
i.e., they coincide with (1.9) up to an elementary factor.

Recall that the Jacobi polynomials are given by 
$$
P_n^{\gamma,\delta}(x)=
\mathrm{const}\cdot\F[-n,n+\gamma+\delta,\delta+1; (1+x)/2] 
$$ 
We observe that
the Romanovski polynomials are analytic continuations of the Jacobi polynomials
with respect to the superscripts.

Many orthogonal systems of this kind are known, see \cite{Ask},
\cite{Les1}--\cite{LW},
\cite{BO}, \cite{Ner-wilson},
see also \cite{Pee}, \cite{BO}, \cite{BO2}
for some applications.  Since these systems are
eigenfunctions of some second order differential or difference operators, we
obtain a good collection of spectral problems (and hence this must give a
collection of new integral transforms).%
\footnote{It seems,that now only difference
operators remain interesting.
 There are 3 types of difference operators related to

a) Shift operator on the lattice $\Bbb Z$:  $Tf(n)= f(n+1)$.

b) Shift operator on the line $\R$: $Tf(x)= f(x+1)$.

c) Shift operator on the line $\R$ in the imaginary direction $Tf(x)=f(x+i).$

It seems, that for the case b), we have infinite multiplicities.  In any case,
variants a), c) are interesting for our purposes.}
 Recently, one of the most
  complicated problems
of this kind
(related to one of families of
hypergeometric ${}_4F_3$-polynomials from \cite{Ner-wilson}) was solved
by Groenevelt \cite{Gro}.

\smallskip

{\bf \punkt Vector-valued bases.}  Fix parameters $0< \alpha<1/2$,
$\beta,p,q\in\R$.  We consider the Hilbert space $Y(\alpha,\beta;q)$, consisting
of $\C^2$-valued functions on half-line $s\ge 0$ with the scalar product
\begin{multline}
 \langle (\phi_1, \phi_2), (\psi_1,\psi_2)\rangle = \\=
\frac{1}{2} \int_0^\infty \begin{pmatrix}\phi_1(s)& \phi_2(s) \end{pmatrix} R(s)
\begin{pmatrix}\ov{\psi_1(s)}\\ \ov{\psi_2(s)} \end{pmatrix} \Bigl|
\frac{\Gamma(1+iq+is)\Gamma(1-iq+is)}{\Gamma(2is)}\Bigr|^2 ds
,\end{multline}
where the matrix $R$ is the same as above (1.6).

Next, define functions $\Xi^{(1)}_n(\alpha,\beta;p,q;s)$,
 $\Xi^{(2)}_n(\alpha,\beta;p,q;s)$ in the variable $s$, given by the formula
 \begin{multline} 
 \Xi^{(1)}_n(\alpha,\beta;p,q;s)=
 \frac{\cos(p-\alpha-iq+i\beta)\pi/2 \,\Gamma(1-\alpha+i\beta)}
 {\Gamma(\alpha-i\beta) \Gamma(1+iq+i\beta) \Gamma(1+iq-\alpha)} \times \\
 \times \FF\Bigr[\begin{matrix} (1-\alpha+i\beta-p+iq)/2-n,\, 1/2+iq+is,
 1/2+iq-is \\ 1+iq+i\beta,1+iq-\alpha \end{matrix};1\Bigr] 
 ;\end{multline}
  $$
 \!\!\!  \!\!\!\!\!\!\!\!\!\!\!\!  \!\!\!  \!\!\!\!\!\!\!\!\!\!\!\!  \!\!\!
 \!\!\!\!\!\!\!\!\!\!\!\!  \!\!\!  \!\!\!\!\!\!\!\!\!\!\!\!  \!\!\!
 \!\!\!\!\!\!\!\!\!\!\!\!  \Xi^{(2)}_n(\alpha,\beta;p,q;s) =
 \Xi^{(1)}_n(-\alpha,-\beta;p,q;x)
  .$$

{\sc Remark.}  The hypergeometric series $\FF(a_1,a_2,a_2;b_1,b_2;1)$ is
 absolutely convergent for $\sum a_i<\sum b_j$ and admits an analytic
 continuation as a meromorphic (sigle-valued) function to arbitrary values of
 the parameters $a_1$, $a_2$, $a_2$, $b_1$, $b_2$.  We understand (1.12) as a
 value of this meromorphic function.  Below (3.4), we represent the function
 $\Xi^{(1)}$ as a hypergeometric series that converges for all interesting for us
 values of the parameters.

\smallskip

{\sc Theorem 1.3.}  {\it The system $(\Xi^{(1)}_n,\Xi^{(2)}_n)$, where $n$
ranges in $\Z$, is an orthonormal basis of the Hilbert space
$Y(\alpha,\beta;q)$.  }

\smallskip

{\sc Remark.}  Emphasis, that the expression (1.12) has the structure
\begin{equation}
 \mathrm{const}\cdot \FF\Bigl[\begin{matrix} -n+h,a+is,a-is\\
a+b, a+c\end{matrix} ;1\Bigr]
\end{equation}
 with
 $$
 \Re a=1/2,\qquad \Re
b=1/2,\qquad c\in\R
.$$
Also, in comparison with  formula (1.4) for
the Hahn
polynomials, the parameter $n$ is shifted by some complex value $h$.


\smallskip

{\sc Remark.}  Similar perturbations exist for other
classical hypergeometric
 systems (Laguerre, Jacobi, Meixner--Pollachek, Meixner),
 they are discussed in
 \cite{Ner-prep}).%
 \footnote{An attempt to find some ways for vector-valued
 special functions starting from spherical functions
 is contained in \cite{GPT}, \cite{Tir},
  see also further
 references in these works. This approach differs from our standpoint
 based on spectral problems with multiple spectra.}

\smallskip

{\bf \punkt Further structure of the paper.}  Theorems 1.1--1.2 are proved in
\S 2, Theorem 1.3 is obtained in \S 3.

In Addendum we discuss perturbations of Laguerre,
Meixner, Meixner--Pollachek, and Jacobi systems.

\bigskip

{\bf\large $\mathbf\S$2.  Spectral decomposition of the hypergeometric differential
operator on the contour $\Re z=1/2$.}

\nopagebreak

\addtocounter{sec}{1} \setcounter{equation}{0} \setcounter{punkt}{0}

\bigskip

A proof of Theorems 1.1--1.2 proposed below
is direct but very tedious.
 We directly apply the Weyl--Titchmarsh--Kodaira theorem
 for an appropriate differential operator.
This theorem with lot of examples is analyzed in Chapter XIII
of  Dunford--Schwartz's book
  \cite{DS}, see also the  Titchmarsh's  book \cite{Tit}.

For the classical index hypergeometric transform
(1.3) several proofs of the inversion
formula are known, see \cite{Koo2}, \cite{Ole},
\cite{VK}, 7.8.8.
It is interesting to find
a shorter proof of Theorem 1.1.

\smallskip

{\bf \punkt Hypergeometric operator.}  Let  $\alpha\ge0$, $\beta\in \R$.
We consider the differential operator
 \begin{equation}
D:=\bigl(\frac14+x^2\bigr)\frac{d^2}{dx^2} + 2x \frac d{dx}+
\frac{(\alpha+i\beta)^2}{4(1/2+ix)}+ \frac{(\alpha-i\beta)^2}{4(1/2-ix)}+\frac14
.\end{equation} 
This operator is formally self-adjoint
in  $L^2(\R)$.  Its resolvent 
 $(D-\lambda)^{-1}$
 is explicitly evaluated below.
 It is well defined for 
  $\Im\lambda\ne 0$. 
  Hence the deficiency indices of $D$
  are  0; this implies the self-adjointness of the
  operator $D$. Our purpose is to construct
  the spectral decomposition
  of the operator  $D$.

Let
 $$r(x):=(\tfrac 12+ix)^{-(\alpha+i\beta)/2} (\tfrac
12-ix)^{-(\alpha-i\beta)/2} .$$
Evaluating directly the differential operator
 $$Bf:=r^{-1}D(rf),$$
 we obtain
  $$ B=
\bigl(\frac14+x^2\bigr)\frac{d^2}{dx^2}+ \bigl(\beta+x(2-2\alpha)\bigr) \frac
d{dx} +\bigl(\alpha-\frac12\bigr)^2 $$
Passing to the complex variable
 $$z=1/2+ix,$$ 
 we obtain the operator
  \begin{equation} A:=-z(1-z)\frac
{d^2}{dz^2}- \bigl(1-\alpha-i \beta- z(2-2\alpha)\bigr) \frac d{dz}-
\bigl(\alpha-\frac12\bigr)^2 .\end{equation} 
Hence the equation
 $Af=\mu^2 f$
transforms to the hypergeometric
equation
 \begin{equation}
\Bigl[z(1-z) \frac{d^2}{dz^2} + \bigl(c-(a+b+1)z)\frac d{dz} -ab\Bigr]\,f=0
\end{equation}
with
 $$ a=\tfrac12-\alpha+\mu,\qquad b=\tfrac12-\alpha-\mu, \qquad
c:=1-\alpha-i\beta $$

Now, we can use the Kummer series
for solutions of the equation
$Df=\mu^2 f$. Below we use the standard notation
 $u_1$,\dots, $u_6$ of
\cite{HTF} for the Kummer solutions
of the standard hypergeometric equation
 (2.3) and use explicit formulae  \cite{HTF}, (2.9.1)--(2.9.24)
 for their expansions in series.

\smallskip

{\bf \punkt Bases in the space of solutions of
the hypergeometric equation.}
We will use four bases in the space of solutions
of the equation
  $(D-\mu^2)f=0$.

\smallskip

{\sc Basis $S_1$, $S_2$.}  Writing the Kummer solutions 
 $u_1$, $u_5$
of the usual hypergeometric equation
 (2.3),  we obtain the following pair of solutions
  $S_1$, $S_2$ 
  of the equation
   (2.1):
     \begin{multline} S_1(\alpha,\beta;
\mu;x):= \bigl(\tfrac 12+ix\bigr)^{-(\alpha+i\beta)/2} \bigl(\tfrac
12-ix\bigr)^{(\alpha+i\beta)/2-\mu-1/2} \times \\ \times \F
\Bigl[\begin{matrix}1/2-\alpha+\mu,1/2-i\beta+\mu\\ 1-\alpha-i\beta
\end{matrix};\, \frac{ix+1/2} {ix-1/2} \Bigr] ;\end{multline} $$
\!\!\!\!\!\!\!\!\!\!\!\!\!\!\!  \!\!\!\!\!\!\!\!\!\!\!\!\!\!\!
\!\!\!\!\!\!\!\!\!\!\!\!\!\!\!  \!\!\!\!\!\!\!\!\!\!\!\!\!\!\!
\!\!\!\!\!\!\!\!\!\!\!\!\!\!\!  \!\!\!\!\!\!\!\!\!\!\!\!\!\!\!
S_2(\alpha,\beta; \mu;x) = S_1(-\alpha,-\beta; \mu;x) .$$
The hypergeometric series
in (2.4) is absolutely convergent
for  $\Im x>0$,
or equivalently, for
$\Re z<1/2$.
On the line  $x\in\R$ the series is conditionally  convergent.
Nevertheless it admits analytical continuation
through the line  $\Im x=0$,
and it is more convenient to think that we consider its analytic continuation.

Further, assume that for  $x=iy$
with  $-1/2<y<1/2$
(or, equivalently, for
$0<z<1$)
 $$
  \bigl(\tfrac12+ix)^\lambda:=e^{\lambda \ln(1/2+ix)} ,\qquad
\bigl(\tfrac12-ix)^\nu:=e^{\nu \ln(1/2-ix)} .
$$
Now we can assume that our solutions $S_1$, $S_2$
are defined
 in the domain $\bigl\{\Re z\le
1/2\bigr\}\setminus[-\infty,0)$.

The solutions
 $S_1$, $S_2$ have asymptotics
  \begin{equation} S_1(x)\sim
\bigl(\tfrac 12+ix\bigr)^{(\alpha+i\beta)/2} ,\qquad S_2(x)\sim \bigl(\tfrac
12+ix\bigr)^{-(\alpha+i\beta)/2} ; \qquad \text{for $x\to i/2$}
.\end{equation}

We also can define  $S_1$ by the formula
 \begin{equation}
 \bigl(\tfrac
12+ix\bigr)^{-(\alpha+i\beta)/2} \bigl(\tfrac 12-ix\bigr)^{-(\alpha-i\beta)/2}
\F \Bigl[\begin{matrix}1/2-\alpha+\mu,1/2-\alpha-\mu\\ 1-\alpha-i\beta
\end{matrix};\, {ix+1/2} \Bigr] \end{equation}
Here the hypergeometric series
converges for  $|1/2-ix|<1/2$, i.e.,
we do not obtain an explicit formula
on the whole line  $x\in\R$.

We mention that
 $$
 S_1(\alpha,\beta;-\mu;x)=S_1(\alpha,\beta;\mu;x)
 .$$

{\sc Basis $T_1$, $T_2$.}  
Writing the pair of the Kummer solutions
 $u_2$, $u_6$,
we obtain the following pair
of solutions
 $T_1$, $T_2$ of the equation
  $(D-\mu^2)f=0$
$$
 T_{1,2}(\alpha,\beta;\mu;x):= S_{1,2}(\alpha,-\beta;\mu;-x) .$$ 
 These solutions 
are defined in the domain 
 $\bigl\{\Re z\ge 1/2\bigr\}\setminus[1,\infty)$.

\smallskip

{\sc Bases  $V_-$, $V_+$  and $W_-$, $W_+$.} 
The Kummer solutions  $u_3$, $u_4$
are defined outside the circular lune $|z|<1$, $|z-1|<1$;
in particular, 
this lune contains the segment 
 $-\sqrt3/2<x<\sqrt3/2$.
Hence  the analytic continuations of
$u_3$, $u_4$ from
the points $z=1/2+i\infty$ and $z=1/2-i\infty$
 to the whole line $\Re z=1/2$
 are different. This gives  the following
 two pairs $V_\pm$, $W_\pm$ of eigenfunctions of (2.1).

 We define the solution
  $V_-$ of the equation
$(D-\mu^2)f=0$
as the solution
that for  $x>\sqrt 3/2$ is given by the formula
\begin{multline}
V_-(\alpha,\beta;\mu;x)= e^{(-1/2+\alpha+i\beta-\mu)\pi i/2}
\bigl(\tfrac 12+ix\bigr)^{-(\alpha+i\beta)/2} \bigl(\tfrac
12-ix\bigr)^{-1/2+(\alpha+i\beta)/2-\mu} \times \\ \times \F\Bigl[\begin{matrix}
1/2-\alpha+\mu,1/2-i\beta+\mu\\ 1+2\mu\end{matrix};\,\frac 1{1/2-ix} \Bigr]
.\end{multline}
Further, assume
 $$ V_+(\alpha,\beta;\mu;x)
=V_-(\alpha,\beta;-\mu;x) .$$
The asymptotics of  $V_{\pm}$
for  $x\to+\infty$ has the form
\begin{equation}
V_{\pm}\sim x^{1/2\pm \mu}, \qquad x\to+\infty .\end{equation}
The solutions $W_\pm$ are determined by the condition
 \begin{equation}
  W_{\pm}\sim
x^{1/2\pm \mu}, \qquad x\to-\infty
.\end{equation}

To obtain a formula for
 $W_-$ for $x<\sqrt 3/2$
 we must change the sign in the argument of the exponential
 function
 in  (2.7). Next,
  $ W_+(\alpha,\beta;\mu;x)
=W_-(\alpha,\beta;-\mu;x) $.

{\bf\punkt Transition matrices.} Define the constants
 \begin{align}
C(\alpha,\beta;\mu)&:= \frac{\Gamma(\alpha+i\beta)\Gamma(1+2\mu)}
{\Gamma(1/2+\alpha+\mu)\Gamma(1/2+i\beta+\mu)} ;\\ \chi(\alpha,\beta;\mu)&:=
e^{-(1/2+\alpha+i\beta-\mu)\pi i/2} 
.\end{align}
In this notations
 (we use  \cite{HTF}, (2.9.37), (2.9.39)),
 \begin{align}
V_-(\alpha,\beta;\mu;x)=& C(\alpha,\beta;\mu) \chi(\alpha,\beta;\mu)
S_1(\alpha,\beta;\mu;x) +\nonumber \\ &+ C(-\alpha,-\beta;\mu)
\chi(-\alpha,-\beta;\mu) S_2(\alpha,\beta;\mu;x) ; \\ V_+(\alpha,\beta;\mu;x)=&
C(\alpha,\beta;-\mu) \chi(\alpha,\beta;-\mu) S_1(\alpha,\beta;\mu;x) +\nonumber
\\ &+ C(-\alpha,-\beta;-\mu) \chi(-\alpha,-\beta;-\mu) S_2(\alpha,\beta;\mu;x)
.\end{align}

Formulae expressing $W_-$, $W_+$
in  $S_1$, $S_2$ can be obtain from 
 (2.12)--(2.13) by the transform
  $\chi\mapsto\chi^{-1}$.

We also need in
an expression of
 $W_\pm$  in terms of  $T_1$, $T_2$.
 Similarly,
 applying  \cite{HTF}, (2.9.38), (2.9.40),
 we obtain
  \begin{align}
W_-(\alpha,\beta;\mu;x)=& -C(\alpha,-\beta;\mu)\chi^{-1}(-\alpha,\beta;-\mu)
T_1(\alpha,\beta;\mu;x) \nonumber -\\&
-C(-\alpha,\beta;\mu)\chi^{-1}(\alpha,-\beta;-\mu) T_2(\alpha,\beta;\mu;x)
\end{align}
To obtain a formula for  $W_+$
we must change sign of  $\mu$
in
$C(\alpha,\beta;\mu)$ and in $\chi(\alpha,\beta;\mu)$.

\smallskip

{\bf\punkt Evaluation of Wronskians.}
Denote by
$\vr(P,Q):=\det\begin{pmatrix}P&Q\\P'&Q'\end{pmatrix}$
the Wronskian of two arbitrary solutions  $P$, $Q$ 
of the equation 
$(D-\mu^2)f=0$.
This expression must have a form 
$$\vr(P,Q)=\frac{\sigma(P,Q)}{1/4+x^2},$$ 
where
 $\sigma(P,Q)$
 is a constant, see
\cite{Kam}, I.17.1.

The value  $\sigma(S_1,S_2)$ 
can be easily evaluated using asymptotics
 (2.5), we obtain 
$$ 
\sigma(S_1,S_2)=i(\alpha+i\beta) 
.$$
 Below, we need in 
$\sigma(V_-,W_-)$.  
The determinant  $\Delta$ of the transition
matrix from the basis 
$(S_1,S_2)$ to the basis  $(V_-,W_-)$ 
can be easily evaluated, this gives
 \begin{equation}
\sigma(V_-,W_-)=\Delta\cdot \sigma(S_1,S_2)= \frac{2\pi i C(\alpha,\beta;\mu)
C(-\alpha,-\beta;\mu)} {\Gamma(\alpha+i\beta) \Gamma(-\alpha-i\beta)} .
\end{equation}

{\bf \punkt Kernel of resolvent.} Now we are ready to write the kernel
$K(x,y;\lambda)$ of
the resolvent
 $$
 R(\lambda):=(D-\lambda)^{-1}
 $$
 of the operator $D$,
  $$
R(\lambda) f(x)=\int_{-\infty}^\infty K(x,y;\lambda) f(y)\,dy. 
 $$ 
 Assume  $\lambda\in\C\setminus(-\infty,0)$.
 Then the solution $V_-$ is an element of 
 $L^2(0,+\infty)$
 and $W_-\in L^2(0,-\infty)$, see  (2.8)--(2.9).
Hence (see  \cite{DS},XIII.3.6), 
$$
 K(x,y;\lambda)=
 \left\{\begin{aligned}
\frac{V_-(x)W_-(y)}{\sigma(V_-,W_-)}, \qquad \text{for $y<x$} \\
\frac{V_-(y)W_-(x)}{\sigma(V_-,W_-)}, \qquad \text{for $x<y$} , \end{aligned}
\right.
$$
the value $\sigma(V_-,W_-)$ was evaluated above  (2.15).

We intend to write an expansion of the differential
operator
 $R$ in eigenfunctions. According prescription
 rising to Weyl's work  \cite{Wey} 
 (see \cite{Tit} 
 and detailed presentation in Dunford, Schwartz
  \cite{DS},
XIII.5.18), 
we must evaluate the jump of the resolvent
on the real line 
\begin{equation} 
\frac 1{2\pi i}\lim_{\epsilon\to 0} \int_{-\infty}^\infty
(R(\lambda+i\epsilon)-R(\lambda-i\epsilon)\bigr)\,d\lambda 
,\end{equation}
and represent it in the form  
\begin{equation} 
Lf(x)= \sum_{i=1,2}\,\sum_{j=1.2}
\int_{-\infty}^\infty \Bigl( \int_{-\infty}^\infty
\sigma_i(x,\lambda)\ov{\sigma_j(y,\lambda))}\,f(y)\, d\mu_{ij}(\lambda)\,
\Bigr)\,dy
 \end{equation} 
for some solutions 
 $\sigma_1$,
$\sigma_2$ of the equation  $Df=\lambda f$. 
Here  $\mu_{ij}(\lambda)$ are (complex-valued)
 measures on  $\R$.  Then  $\mu_{ij}$ is the spectral measure.
 Precisely, the operator
 $$ 
 f\mapsto \Bigl( \int_{-\infty}^\infty f(x)\ov{\sigma_1(x)}\,dx,
\int_{-\infty}^\infty f(x)\ov{\sigma_2(x)}\,dx \Bigr)
 $$
 is a unitary operator from
  $L^2(\R)$ to the space of  $\C^2$-valued functions
  $(\phi_1(\lambda),\phi_2(\lambda))$
  with the scalar product 
   $$ 
   \bigl\langle
(\phi_1(\lambda),\phi_1(\lambda), (\psi_1(\lambda),\psi_1(\lambda) \bigr\rangle=
\sum_{i=1,2}\,\sum_{i=1,2}\int_{-\infty}^\infty \phi_i(\lambda)\ov{
\phi_j(\lambda) } \,d\mu_{ij}(\lambda)
 $$ 
 (in particular, the right-hand side
 of this equality is a positive definite
 scalar product).

\smallskip

{\bf \punkt Formula for resolvent.} 
To
evaluate the jump of the resolvent for
 $\lambda\ge0$ and for  $\lambda\le0$, we are need
 of two
 explicit expressions for the kernel
  $K(x,y;\lambda)$ of the resolvent.

Expressing $V_-$ via  $S_1$, $S_2$
and $W_-$ via  $T_1$, $T_2$, we obtain
the following expression for the resolvent
 \begin{multline}
 K(x,y;\lambda)= \frac
{-1}{2\pi } \Gamma(\alpha+i\beta) \Gamma(-\alpha-i\beta) \biggl[
\frac{C(\alpha,-\beta;\sqrt \lambda)}{C(-\alpha,-\beta;\sqrt \lambda)}
e^{(-1/2+\alpha-\sqrt \lambda)\pi i} S_1(x) T_1(y) 
+ \\+
\frac{C(-\alpha,\beta;\sqrt \lambda)}{C(-\alpha,-\beta;\sqrt \lambda)}
e^{(-1/2+i\beta-\sqrt \lambda)\pi i} S_1(x) T_2(y) + \frac{C(\alpha,-\beta;\sqrt
\lambda)}{C(\alpha,\beta;\sqrt \lambda)} e^{(-1/2-i\beta-\sqrt \lambda)\pi i}
S_2(x) T_1(y) +\\+ \frac{C(-\alpha,\beta;\sqrt \lambda)}{C(\alpha,\beta;\sqrt
\lambda)} e^{(-1/2-\alpha-\sqrt \lambda)\pi i} S_2(x) T_2(y) \Bigr]
,\end{multline} 
where $C(\dots)$ is defined by  (2.10).

Expressing  $V_-$, $W_-$ by  $S_1$, $S_2$, 
we obtain
 \begin{multline}
K(x,y;\lambda)= \frac{\Gamma(\alpha+i\beta) \Gamma(-\alpha-i\beta)} {2\pi
C(\alpha,\beta;\sqrt\lambda)\,C(-\alpha,-\beta;\sqrt\lambda)} \times \\ \times
\Bigl[ C(\alpha,\beta;\sqrt\lambda))\, \chi(\alpha,\beta;\sqrt\lambda) \, S_1(x)
+ C(-\alpha,-\beta;\sqrt\lambda))\, \chi(-\alpha,-\beta;\sqrt\lambda)) \, S_2(x)
\Bigr] \times \\ \times \Bigl[ C(\alpha,\beta;\sqrt\lambda))\,
\chi^{-1}(\alpha,\beta;\sqrt\lambda) \, S_1(x) +
C(-\alpha,-\beta;\sqrt\lambda))\, \chi^{-1}(-\alpha,-\beta;\sqrt\lambda)) \,
S_2(x) \Bigr]
 ,\end{multline} 
 where $\chi$ is given by  (2.11).

These two expressions are meromorphic in the plane
 $\lambda\in\C$,
with a cut on the negative semiaxis  $\lambda<0$.

\smallskip

{\bf \punkt Jump of resolvent for $\lambda>0$.} 
Evaluation of 
\begin{equation} 
L^{[0,\infty)}:= \frac 1{2\pi i}\lim_{\epsilon\to 0}
\int_{0}^\infty \bigl(R(\lambda+i\epsilon)-R(\lambda-i\epsilon)\bigr)\,d\lambda
\end{equation}
is reduced to an evaluation of the residues.
Expression  (2.19) has simple poles 
at the points 
 $\lambda=(\alpha-k-1/2)^2$ 
 for integer  $k$,
satisfying the condition
 $0\le k\le \alpha-1/2$.  In addition, 
only the coefficient at  $S_1(x)S_1(y)$
has nonzero residues, three other coefficients 
are holomorphic at these points.

Thus the integral operator
 (2.20) equals 
  \begin{multline*}
L^{[0,\infty)} f(x)= \frac 1{2\pi} \Gamma(\alpha+i\beta) \Gamma(\alpha-i\beta)
\sum_{0\le k <\alpha-1/2} \frac{(2\alpha-2k-1)\,(1-\alpha-i\beta)_k}
{\Gamma(2\alpha-k)\,(1-\alpha+i\beta)_k\,k!}  \times \\ \times
\int_{-\infty}^\infty S_1(\alpha,\beta;\alpha-k-1/2;x)
S_1(\alpha,\beta;\alpha-k-1/2;y) f(y)\,dy
 \end{multline*}
 Writing  $S_1$  at the form  (2.6), 
 we obtain its representation in the terms of the Romanovski
 polynomials.
 Applying identity  \cite{HTF}, (10.8.16)
 for Jacobi polynomials, we convert our expression
 to the form 
 \begin{multline*} 
 \frac 1{2\pi}
\Gamma(\alpha+i\beta) \Gamma(\alpha-i\beta) \sum_{0\le k <\alpha-1/2}
\frac{2\alpha-2k-1} {\Gamma(2\alpha-k)\,k!}  S_1(\alpha,\beta;\alpha-k-1/2;x)
\times\\ \times \int_{-\infty}^\infty \ov {S_1(\alpha,\beta;\alpha-k-1/2;y)}
f(y)\,dy
 \end{multline*}
and this gives the required expression.

\smallskip

{\bf \punkt Jump of resolvent for  $\lambda<0$.} 
We must evaluate 
$$
K(x,y;\sqrt\lambda)-K(x,y;-\sqrt\lambda),
$$
 assuming that  $\lambda$ is negative real
 and  $\Im \sqrt\lambda>0$.

For definiteness, we present the calculation
of the coefficient at 
 $S_1(x)T_1(y)$.
In three remaining cases,
the evaluation
is almost identical%
 \footnote{This is explained by natural
 symmetries
  $(\alpha,\beta)\sim (-\alpha,-\beta)
\sim (i\beta, -i\alpha)$ of the equation  (2.1).}.

Thus,  (see  (2.18)), we intend to find
 \begin{multline*} 
 \frac1{-2\pi}
\Gamma(\alpha+i\beta)\Gamma(-\alpha-i\beta) \times \\ \times \Bigl\{
\frac{C(\alpha,-\beta;\sqrt\lambda)} {C(-\alpha,-\beta;\sqrt\lambda)}
e^{-(1/2-\alpha+\sqrt\lambda)\pi i} - \frac{C(\alpha,-\beta;-\sqrt\lambda)}
{C(-\alpha,-\beta;-\sqrt\lambda)} e^{-(1/2-\alpha-\sqrt\lambda)\pi i} \Bigr\}
.\end{multline*}
After direct cancellations, we obtain 
\begin{multline*}
 \frac1{-2\pi} \Gamma(\alpha+i\beta)\Gamma(\alpha-i\beta)
\times\\ \times \Bigl\{\frac{\Gamma(1/2-\alpha+\sqrt\lambda)}
{\Gamma(1/2+\alpha+\sqrt\lambda)} e^{-(1/2-\alpha+\sqrt\lambda)\pi i} +
\frac{\Gamma(1/2-\alpha-\sqrt\lambda)} {\Gamma(1/2+\alpha-\sqrt\lambda)}
e^{-(1/2-\alpha-\sqrt\lambda)\pi i} \Bigr\} 
.\end{multline*}
We convert this expression to
\begin{multline*} \frac1{-2\pi}
\Gamma(\alpha+i\beta)\Gamma(\alpha-i\beta) \Gamma(1/2-\alpha+\sqrt\lambda)
\Gamma(1/2-\alpha-\sqrt\lambda) \times\\ \times
\Bigl\{\frac{e^{-(1/2-\alpha+\sqrt\lambda)\pi i}}
{\Gamma(1/2+\alpha+\sqrt\lambda)\Gamma(1/2-\alpha-\sqrt\lambda)} -
\frac{e^{-(1/2-\alpha-\sqrt\lambda)\pi i}}
{\Gamma(1/2+\alpha-\sqrt\lambda)\Gamma(1/2-\alpha+\sqrt\lambda)} \Bigr\}
.\end{multline*}
Now, we will transform only the expression
in the curly brackets
 $$ 
 \pi\Bigl\{\dots\Bigr\}= \cos\bigl(\pi (\alpha+\sqrt\lambda)\bigr)
e^{-(1/2-\alpha+\sqrt\lambda)\pi i} - \cos\bigl(\pi (\alpha-\sqrt\lambda)\bigr)
e^{-(1/2-\alpha-\sqrt\lambda)\pi i}
 .$$
 Next, we apply the Euler formulae for 
$\cos$, collect similar terms in the linear combination
of exponential functions,
and again apply the Euler formulae. 
We obtain 
 $$
 i\sin (2\sqrt\lambda\pi i)
 ,$$
 and this finishes 
the evaluation of the coefficient at  $S_1(x)T_1(y)$.

\smallskip

After this, for the jump of the resolvent
 $L^{(-\infty,0]}$  on the semiaxis 
  $(-\infty,0]$,
  we obtain the expression of the form 
   $$
    L^{(-\infty,0)}f(x)= \sum_{i=1,2}
\sum_{j=1,2} \int_{-\infty}^0 d\lambda\biggl\{ \theta_{ij}(\lambda)
\int_{-\infty} S_i(x)T_j(y)f(y)dy\biggr\}
 $$
 where  $\theta_{ij}$ are some explicit products of 
  $\Gamma$-functions.

Further,  observe
that for  $\lambda<0$ and  $y\in\R$
we have 
$T_j(y)=\ov{S_j(y)}$. 
We obtain an expression of the form  (2.17), 
and this finishes the evaluation of the spectral decomposition
of the operator
 (2.1).

The parameter  $s$ from 1.4 is  $\sqrt\lambda$.

\bigskip

{\bf\large $\mathbf\S$ 3.  Construction of bases}

\bigskip

\addtocounter{sec}{1} 
\setcounter{equation}{0} 
\setcounter{punkt}{0}

{\bf \punkt One basis in $L^2(\R)$.}  
Fix real parameters  $p$ and
$q$.

\smallskip

{\sc Lemma 3.1.}  {\it The system of functions
\begin{equation} 
r^{(n)}_{p,q}(x):=
\bigl(\tfrac12 +ix\bigr)^{-1/2-n-(p+iq)/2} \bigl(\tfrac12
-ix\bigr)^{-1/2+n+(p-iq)/2} 
,\end{equation}
where  $n$ ranges in $\Z$,
is an orthonormal basis of $L^2(\R, dx/2\pi)$.}

{\sc Proof.}  Pass to a new variable
 $\psi\in[0,2\pi]$ by the formula
$$
e^{i\psi}=\frac{1/2 +ix} {1/2 -ix},\qquad d\psi=\frac{dx}{1/4+x^2}
.$$
Our system of functions transforms to
$$
 e^{-in\psi} \cdot e^{-ip\psi/2} (2\cos
\psi/2)^{iq}
.$$
We obtain the standard orthogonal system
 $e^{-in\psi}$
up to multiplication by
a function whose absolute value is 1.
 \hfill$\square$

Emphasis also that orthogonality of our system can be obtained
directly (since the scalar products can be easily
evaluated using residue calculus).

{\bf \punkt Construction of bases.} The bases (1.13)
can be easily obtained, if we apply the double index hypergeometric
transform to the orthogonal system  (3.1).

Let us explain, how to perform the calculation.
We must find
 \begin{equation}
\int_{-\infty}^\infty r^{(n)}_{p,q}(x) \ov{Q_1(\alpha,\beta;x;s)}\,dx
\end{equation} 
where the function  $Q_1$ is defined by the formula  (1.6).
For this purpose, we expand  $\F$ from formula  (1.6)
in the hypergeometric series in powers of $$
\frac{ix+1/2}{ix-1/2} 
$$
Then we integrate it termwise using the Cauchy beta-integral
(see  \cite{PBM}, v.1, (2.2.6.31)
or \cite{AAR}, Chapter 1, ex.  13)
\begin{equation}
\int\nolimits_{-\infty}^\infty 
\frac {dx} {(1/2+ix)^\sigma
(1/2-ix)^\tau} = \frac{2\pi\Gamma(\sigma+\tau-1)}{\Gamma(\sigma)\Gamma(\tau)}
\end{equation}

As a result, we obtain the following expression for 
 (3.3) 
\begin{multline}
 \frac{2\pi \Gamma(1/2+iq-is)}
{\Gamma(\alpha-i\beta)\Gamma((p+iq-\alpha+i\beta)/2+n+1-is)
\Gamma((-p+iq+\alpha-i\beta+1)/2-n)} \times \\ \times
\FF\Bigl[\begin{matrix}1/2-\alpha-is, 1/2+i\beta-is,
1/2+(p-iq-\alpha+i\beta)/2+n\\ 1-\alpha+i\beta, (p+iq-\alpha+i\beta)/2+n+1-is
\end{matrix} ; 1\Bigr] 
\end{multline} 
This expression is a variant of the final answer,
its more symmetric form  (1.13)
can be obtained by Thomae transformation
see  \cite{AAR}, Corollary  3.3.6, \cite{PBM},
v.3, (7.4.4.2).
$$ 
\FF\Bigl[\begin{matrix} a,b,c\\d,e\end{matrix};1\Bigr]=
\frac{\Gamma(d)\Gamma(e)\Gamma(r)} {\Gamma(a)\Gamma(b+r)\Gamma(c+r)}
\FF\Bigl[\begin{matrix} d-a,e-a,r\\b+r,c+r\end{matrix};1\Bigr] 
,$$ 
where 
$r=d+e-a-b-c$.




{\sf Math.Physics group, 
 Institute of Theoretical and Experimental Physics, %
\linebreak
B.Cheremushkinskaya, 25, Moscow 117 259, Russia}


{\tt neretin@mccme.ru

\&

ESI, Wien (November--December, 2003)
}




\newcommand{\tha}{\mathop{\rm th}\nolimits}

\def\cL{{\EuScript L}}
\def\cM{{\EuScript M}}
\def\cP{{\EuScript P}}
\def\cS{{\EuScript S}}

\renewcommand\FF{\,\,{}_2F_1}






\newpage

\begin{center}
\Large\bf

Addendum.

Perturbations of some

classical hypergeometric orthogonal systems

\end{center}

Here we discuss orthogonal systems
that can be obtained from
some classical orthogonal system
$T_n$ by transformation $n\mapsto n+\theta$,
there $\theta\in\C$ is fixed.

The  most of facts formulated
in this addendum are known
or semi-known.

\smallskip

{\bf A.1. Spectral problems for one-parametric subgroups in
$\SL(2,\R)$.}
Consider a unitary irreducible representation
of the universal covering of the group $\SL(2,\R)$.
There are three (up to a conjugation) one-parameter subgroups
in $\SL(2,\R)$, namely

 -- an elliptic subgroup $K=\SO(2)$,
consisting of orthogonal matrices;

   -- a hyperbolic subgroup $D$ consisting
   of diagonal matrices $\begin{pmatrix} a&0\\0&a^{-1}\end{pmatrix}$, where $a>0$

   -- a parabolic subgroup $N$ consisting of matrices
   $\begin{pmatrix} 1 &t\\0&1\end{pmatrix}$

 Consider the spectral decomposition
 of our representation with respect to
  one of these subgroups
 (i.e., consider a realization of our representation in some space $L^2$,
 where our subgroup acts by multiplication by functions).
 Operators of representation of $\SL(2,\R)$ in these models
 were obtained in the books of Vilenkin
\cite{AVil}, Chapter 7,
and of Vilenkin and Klimyk \cite{AVK},
Sections 6.8, 7.2, 7.6%
\footnote{To be precise, they usually consider   the group
$\SL(2,\R)$ itself,
extension of their formulae to
the universal covering group is not a difficult problem.}
\footnote{Likely, the earliest references that can be attributed
to this subject are \cite{AKep}, \cite{AML}.}

Also  consider some subgroup in $\SL(2,\R)$
 conjugated to $\SO(2)$, consider its eigenbasis
 in the space of representation, and
 consider the image of this eigenbasis in our space
 $L^2$. Obviously, we will obtain an orthogonal basis
 in the space $L^2$.
 Considering the spectral decomposition
 of a highest weight representation
 with respect to $K$, $D$ or $N$,
 we obtain respectively Meixner polynomials, Meixner-Pollachek
 polynomials and Laguerre polynomials
 (see \cite{AVK}, 6.8.3, 7.7.11, 7.7.8).

Starting from a representation of a principal or complementary series,
we obtain some nonpolynomial bases. These cases are discussed below in A.2-A.4.

\smallskip

{\sc Remark.}
In A.2 -- A.4, we omit representations from
our discussion and formally use
only elementary analytic tools.
But the representation theory allows to
propose simple and effective actions. In fact,
initial orthogonal bases,
that we use below, are eigenbases for elliptic subgroups
(for representations realized in functions
on line or circle, see, for instance, \cite{AVK}),
the Fourier transform is the spectral expansion
with respect to the parabolic subgroup $N$,
the Mellin transform is the spectral expansion
with respect to hyperbolic subgroup $D$, and
expansion in Laurent series corresponds
to expansion in eigenvectors of $K$.

\smallskip

{\bf A.2. Perturbation of Laguerre polynomials.}
Recall that the Kummer's confluent hypergeometric function
 ${}_1F_1$ is given by
$$
{}_1F_1[a;c;x]=
\sum_{n=0}^\infty \frac{(a)_n}{(c)_n n! } x^n;\\
$$
The confluent hypergeometric  function
of Tricomi (for details, see \cite{AHTF1},6.5.7
or \cite{ASla1}) is given
$$
\Psi(a,c;x)=\frac{\Gamma(1-c)}{\Gamma(a-c+1)} {}_1 F_1[a;c;x] +
   \frac{ \Gamma(c-1)}{\Gamma(a)}x^{1-c}{}_1 F_1[a-c+1;2-c;x]
$$
%

\smallskip

Fix $\alpha\in\C$. Consider the orthogonal basis
$$
\phi_n(\alpha;x)=\bigl(1/2+ix\bigr)^{-n-1/2-\alpha}
            \bigl(1/2-ix\bigr)^{n-1/2+\ov\alpha}
$$
in $L^2(\R)$, where $n$ ranges in $\Z$.

Considering the Fourier transforms of these functions, we
obtain (see \cite{AIT1}, III.3.2.12; see similar integrals
in \cite{ASla1})
 the following proposition.

 \smallskip

{\sc Proposition A.1.}
{\it The functions
 $\cL^{\alpha,\ov\alpha}_n(y)$ on $\R$
given by
$$
\cL^{\alpha,\ov\alpha}_n(y)=
\left\{
\begin{aligned}
\frac 1{\Gamma(-n+1/2-\ov\alpha)}
  \Psi(1/2+n+\ov\alpha,1+\ov\alpha-\alpha ;y),\quad \text{for $y>0$}
\\
\frac 1{\Gamma(n+1/2+\alpha)}
  \Psi(1/2-n-\alpha,1+\ov\alpha-\alpha ;-y),\quad \text{for $y<0$}
\end{aligned}
\right.
$$
form an orthonormal basis in $L^2$
on $\R$ with respect to the weight
 $\exp(-|y|)$.}

\smallskip

Now, consider $\alpha$, $\beta\in\R$, let
$$|\alpha-\beta|<1/2, \qquad
\text{and $\alpha,\beta\ne\pm 1/2$.}
$$
Define a weight $\mu_{\alpha,\beta}(y)\,dy$  on the line $\R$
by the formula
$$
\mu_{\alpha,\beta}(y)
=
|y|^{\beta-\alpha} e^{-y}\times
\left\{
\begin{matrix} \cos(\pi\beta),\quad\text{for $y>0$}\\
               \cos(\pi\alpha),\quad\text{for $y<0$}
\end{matrix}
\right.
.$$

Define functions    $\cL^{\alpha,\beta}_n$ on $\R$
by
$$
\cL^{\alpha,\beta}_n(y)
=\left\{
\begin{aligned}
\frac 1{\Gamma(1/2-n-\alpha)}
    \Psi(1/2+n+\beta, 1-\alpha+\beta;y)
          ,\quad\text{for $y>0$}   \\
\frac 1{\Gamma(1/2+n+\beta)}
    \Psi(1/2-n-\alpha, 1-\alpha+\beta;y)
                ,\quad\text{for $y<0$}
\end{aligned}
.\right.
$$

{\sc Proposition A.2.} {\it
 The functions $\cL^{\alpha,\beta}_n$,
 where $n$ range in $\Z$, form an orthogonal basis
 in the space  $L^2(\R, \mu(y)\,dy)$. Moreover}
$$
\int_{-\infty}^\infty
\bigl|\cL^{\alpha,\beta}_n(y)\bigr|^2 \mu(y)\,dy=
\frac{(-1)^n\pi}
{\Gamma(1/2+n+\beta)\Gamma(1/2-n-\alpha)}
$$

To prove orthogonality, we write explicitly
the orthogonality relations
for $\cL^{\alpha,\ov\alpha}_n(y)$
and apply the Kummer formula
$$
\Psi(a,c;x)=x^{1-c}\Psi(a-c+1,2-c;x)
$$
After this, we write the analytic continuation
of orthogonality relations with respect to $\alpha$, $\ov\alpha$.

Proposition A.1 corresponds to the principal series of representations
of $\SL(2,\R)$
and Proposition A.2 corresponds to the complementary series.
The  cases $\alpha=1/2$ and $\beta=1/2$
correspond to highest weight and lowest weight representations.
In these cases, the system $\cL^{\alpha,\beta}_n$
degenerates to the usual Laguerre (Sonin) system on half-line.

On these functions $\cL_n$, see  Groenevelt \cite{AGro}.

\smallskip

{\bf A.3. Perturbation of Meixner polynomials.}
We use the following nonstandard notation
for the Gauss hypergeometric function.
$$
\FF^*[a,b;c;z]:=
\frac{\Gamma(a)\Gamma(b)}{\Gamma(c)}
\FF[a,b;c;z]    =
\sum_{k=0}^\infty\frac{\Gamma(a+k)\Gamma(b+k)}{\Gamma(c+k)\,k!} z^k
$$
If $c=-l=0,-1,-2,\dots$, when $\FF[a,b;c;z]$
has poles,
but
\begin{multline*}
\FF^*[a,b;c;z] =
\sum_{m=0}^\infty\frac{\Gamma(a+l+m+1)\Gamma(b+l+m+1)}{m! \,(m+l)!} z^k
=         \\  =
\frac{\Gamma(a+l+1)\Gamma(b+l+1)}{l!} \FF[a+l+1,b+l+1;l+1;z]
\end{multline*}
(first $(l+1)$ summands in the series for $\FF$ are zero).

Consider the following orthogonal system
$$
h_n(z):=z^n (\ch t -z\sh t)^{-1/2-\alpha-n}
            (\ch t- z^{-1}\sh t)^{-1/2+\alpha+n}
$$
in $L^2$ on the circle $|z|=1$; here $n$ ranges in $\Z$.
Expanding these functions in Fourier (Laurent) series,
we obtain the following orthogonality relations
 for their Fourier coefficients.

 \smallskip

{\sc Proposition A.3.} {\it
For  fixed $\alpha\in\C$, $t\in\R$, the functions
$$
\cM_n^{\alpha,\ov\alpha}(j)=
\frac {(\tha t)^{j-n}}{\Gamma(1/2-\ov \alpha-n)}
\FF^*\Bigl[\begin{matrix} 1/2+\alpha+j, 1/2-\ov\alpha-n\\
                           j-n+1\end{matrix}
                         ;\,\tha^2 t
      \Bigr]
$$
where $n$ ranges in $\Z$,
form an orthogonal basis in
 $l^2(\Z)$.
The inner squares of all the elements of
the basis are  $\ch^{-2} t$.}

\smallskip

{\sc Proposition A.4.} {\it
Fix $t\in \R$, fix $\alpha$, $\beta\in\R$
such that
$$|\alpha|<1/2,\quad |\beta|<1/2$$
Consider the space $H_{\alpha,\beta}$
 of two-side sequences  $x(j)$
with the scalar product
$$
\langle x, y\rangle=
\sum_{j=-\infty}^\infty
  x(j)\ov {y(j)} \cdot\frac{\Gamma(1/2+\alpha+j)}
                         {\Gamma(1/2+\beta+j)}
$$
Then the system of functions
$$
\cM^{\alpha,\beta}_n(j):=
\frac{(\tha t)^{j-n}}
     {\Gamma(1/2-\beta-n)} \FF^*\Bigl[
      \begin{matrix}1/2+\alpha+j,1/2-\alpha-n\\
                       j-n-1
      \end{matrix}; \,\tha^2 t\Bigr]
$$
forms an orthogonal basis in  $H_{\alpha,\beta}$,
and}
$$
\sum_{-\infty}^\infty
\bigl|\cM^{\alpha,\beta}_n(j) \bigr|^2
  \cdot\frac{\Gamma(1/2+\alpha+j)}
          {\Gamma(1/2+\beta+j)}
 =\ch^{-2} t \,\Gamma^2(1/2-\alpha-n)
$$

 {\sc Proof.} Consider the following
 biorthogonal system on the circle
\begin{align*}
u_n(z)=z^n (\ch t -z\sh t)^{-1/2-\alpha-n}
            (\ch t- z^{-1}\sh t)^{-1/2+\beta+n};
\\
v_n(z)=z^n (\ch t -z\sh t)^{-1/2-\ov\beta-n}
          (\ch t- z^{-1}\sh t)^{-1/2+\ov\alpha+n}
\end{align*}
Expanding these function in Fourier
(or Laurent) series, we obtain
a biorthogonal system of two-side sequences.
Images $\widehat u_n$, $\widehat v_n$ of
$u_n$, $v_n$ seem different at the first
glance.
But applying the  identity (see \cite{AHTF1}, (2.1.23))
$$
\FF(a,b;c;z)=(1-z)^{c-a-b}\FF(c-a,c-b;c;z)
$$
we observe that $\widehat u_n$,
 $\widehat v_n$ coincide modulo
a simple elementary factor.
Thus we obtain orthogonality relations
for $u_n$.
\hfill $\square$

\smallskip

These facts are contained in Vilenkin--Klimyk \cite{AVK}
(they discuss only two parametric family of bases),
apparently, these bases firstly
appeared in \cite{AVKA}.

\smallskip

{\bf A.4. Perturbation of
the Meixner--Pollachek system.}
Consider the Mellin transform in $L^2(\R)$, i.e.,
for $g\in L^2(\R)$, we define two functions
\begin{align*}
 \widehat g_1(s)&:=\int_0^\infty x^{is-1/2} g(x)\,dx    \\
 \widehat  g_2(s)&:=\int_{-\infty}^0 (-x)^{is-1/2} g(x)\,dx
\end{align*}
on $\R$.
We have the following condition of unitarity
$$
\int_{-\infty}^\infty
\bigl|g(x)\bigr|^2 \,dx=
\frac 1{2\pi} \Bigl\{
\int_{-\infty}^\infty \bigl| \widehat g_1(s)\bigr|^2 dx
+
\int_{-\infty}^\infty \bigl| \widehat g_2(s)\bigr|^2 dx
\Bigr\}
$$

Fix $\alpha=\sigma+i\tau\in\C$ and $\phi\in(0,\pi)$.
Consider the   orthogonal
 basis in $L^2(\R)$ given by
 $$h_n(x)=(1+x e^{i\phi})^{-1/2-\sigma-i\tau-n}
        (1+x e^{-i\phi})^{-1/2+\sigma-i\tau+n}
        x^{i\tau}
$$
We evaluate Fourier transform of these
functions (see \cite{APBM1}, 2.2.6.24) and obtain
the following construction.

Consider the Hilbert space $V_\alpha$,
whose elements are pairs of functions  $(f_1,f_2)$,
depending in the variable $s\in\R$,
and the inner product is
\begin{multline*}
\bigl\langle (f_1,f_2),\, (g_1,g_2)\bigr \rangle=
\frac{|\Gamma(1/2+i\tau+is) \Gamma(1/2+i\tau-is)|^2}
     {2\sin\phi|\Gamma(1+2i\tau)|^2}
 e^{(s+\tau)\phi}
\times \\ \times
\Bigl\{
\int_{-\infty}^\infty f_1(s)\ov{g_1(s)}\,ds
+ e^{-\pi s}\int_{-\infty}^\infty f_2(s)\ov{g_2(s)}\,ds
\Bigr\}
\end{multline*}

Consider the hypergeometric series
$$ \FF\Bigl[\begin{matrix}
       1/2+is+i\tau,1/2-\sigma+i\tau-n\\
        1+2i\tau \end{matrix};\, u
\Bigr]
$$
We consider values of analytic continuation of this series
at the point
$$u=1-e^{-2i\phi}, \qquad\text{where $0<\phi<\pi$}$$
By $G_n^1(s,\phi)$ we denote the value
obtained after passing the way $h(\theta):=1-e^{-2i\theta}$,
 where $\theta \in [0,\phi]$.
By $G_n^2(s;\phi)$ we obtain the result
of passing   $h(\theta):=1-e^{2i\theta}$,
where $\theta \in [0,\pi-\phi]$.

                       \smallskip

{\sc Proposition A.5}
{\it  The system $(G_n^1(s), G_n^2(s) )$
is an orthonormal basis of the space
 $V_\alpha$.}

\smallskip

Thus, we obtain a $\C^2$-valued orthogonal system
(as above in Section 3).

\smallskip

{\bf A.5. On properties of perturbed systems.}
Standard properties of Laguerre, Meixner-Pollachek,
and Meixner polynomials (see \cite{AHTF2} or \cite{AKS})
can be easily extracted from the theory of highest weight
representations of $\SL(2,\R)$; this is explained
in Vilenkin--Klimyk \cite{AVK}.
It is clear that their methods remain
valid for the perturbed
systems.

One of nice facts that can be obtained
in this way is surviving of the Meixner generating
function (see \cite{AMei}).
For instance, consider the
 Poisson kernel for perturbed
Meixner--Pollachek system,
$$
{\mathcal K}_{\alpha,\phi}(s,t;\theta)=
\sum_{n=-\infty}^\infty
\begin{pmatrix}
G_n^1(s)\,\ov{G_n^1(t)}\,&
G_n^1(s)\,\ov{G_n^2(t)}\\
G_n^2(s)\,\ov{G_n^1(t)}\,&
G_n^2(s)\,\ov{G_n^2(t)}
\end{pmatrix}
\, e^{in\theta}
$$
This sum admits an explicit evaluation, since
${\mathcal K}(s,t;\theta)$ are the kernels of
integral operators corresponding
to matrices
$$g(\theta)=
\begin{pmatrix} \cos\theta&\sin\theta\\
                 -\sin\theta&\cos\theta
		 \end{pmatrix}
\in \SL(2,\R)
$$
in the $D$-diagonal model of the representation.
The latter kernels (for $2\alpha\in \Z$) were evaluated by
Vilenkin \cite{AVil}, (7.4.1.13)--(7.4.1.16)
(see also \cite{AVK}, (7.2.1.12)--(7.2.1.15)).

The same is valid for the perturbed Meixner system
and perturbed Laguerre system
(for the Laguerre system, this gives an extension
of
 Myller-Lebedeff  bilinear generating
function, see \cite{AHTF2}, \cite{AML}).

\smallskip

{\bf A.5. Perturbed Jacobi systems
and Molchanov's singular boundary problems.}
Let $-1<\alpha<1$, $\beta>-1$.
Fix $\theta$, such that $0<\theta<1$.
We also assume $0<\theta+\alpha<1$.

{\it Let $n-\theta$ be integer.}
We consider a function $\Phi_n$ on the half-line $x>0$
given by
$$
\Phi_n(x)=
\left\{
\begin{aligned}&
\FF\Bigl[\begin{matrix} -n,n+\alpha+\beta+1\\ \beta+1 \end{matrix};x\Bigr],
    \qquad\qquad \qquad\qquad\text{ for $0<x<1$}
    \\
& x^{-n-\beta-1}(x-1)^{-\alpha}
 \Gamma\Bigl[
 \begin{matrix}
 n+\alpha+\beta+1,\,\alpha+n+1\\
 2n+\alpha+\beta+2
 \end{matrix}
 \Bigr]\times \\
& \qquad\qquad\times
\FF\Bigl[\begin{matrix} n+\alpha+\beta+1,2n+\alpha+1
    \\ 2n+\alpha+\beta+2\end{matrix};\frac 1x\Bigr],
     \qquad \text{ for $x>1$}
\end{aligned}
\right.
$$

{\sc Remark.} These two functions are Kummer solutions
\cite{AHTF1}, (2.9.11), (2.9.14) of the same
differential hypergeometric
equation.

\smallskip

{\sc Theorem A.6.} {\it The system
$\Phi_n$, where $n-\theta\in\Z$ and $2n+\alpha+\beta+1>0$,
  is orthogonal
with respect to the weight}
$$
w(x)=
\left\{
\begin{aligned}&
x^\beta(1-x)^\alpha,\qquad\qquad\qquad\qquad \text{for $0<x<1$};
\\ &
\frac{\sin(\pi\theta)}
     {\sin\pi(\alpha+\theta)}
x^\beta(x-1)^\alpha\qquad \qquad\text{for $x>1$}
\end{aligned}
\right.
$$
{\it and}
$$
\int_0^\infty |\Phi_n(x)|^2 \,dx=
\frac{\Gamma^2(\beta+1)\Gamma(1+n)\Gamma(1+n+\alpha)} {(\alpha+\beta+2n+1)\Gamma(\beta+1+n)\Gamma(n+\alpha+\beta+1)}
$$

{\sc Remark.} {\it We emphasis that the orthogonal system $\Phi_n$
is not a basis.}

\smallskip

{\sc Proof.}
Denote by $H(x)$  the Heaviside
one-step function,
$H(x)=1$ for $x>0$,
and $H(x)=0$ for $x<0$.
Consider two
Mellin--Barnes integrals
\begin{multline*}
K_1(x)=\frac 1{2\pi i}
\Gamma \begin{bmatrix} \beta+1,\, 1+m+\alpha\\
                            \beta+m+1 \end{bmatrix}
    \int_{-i\infty}^{+i\infty}
    \Gamma\begin{bmatrix} s,\,\beta+1+m-s\\
                          m+\alpha+1+s, \beta+1-s
                       \end{bmatrix}
     x^{-s} ds
     =\\ =
\FF\Bigl[ \begin{matrix}
  \beta+1+m,\, -m-\alpha\\ \beta+1
   \end{matrix} ; x\Bigr]\, H(1-x)
   +\\
   + x^{-\beta-1-m}
   \Gamma\begin{bmatrix} \beta+1,\, 1+m+\alpha\\
                        -m, 2m+\alpha+\beta+2
          \end{bmatrix}
  \FF\Bigl[ \begin{matrix} \beta+1+m,\, m+1\\
                            2m+\alpha+\beta+2\end{matrix}
                       ;\frac 1x \Bigr] H(x-1)
\end{multline*}

\begin{multline*}
K_2(x)=\frac 1{2\pi i}
\Gamma
  \begin{bmatrix}2n+\alpha+\beta+2,\, -\alpha-n\\
  n+\alpha+\beta+1
  \end{bmatrix}
\int_{-i\infty}^{+i\infty}
\Gamma\begin{bmatrix}
\alpha+n+s,\, \beta+1-s\\
s,\, n+\beta+2-s
\end{bmatrix}  \,  x^{-s} ds
=\\=
x^{n+\alpha} \FF\Bigl[ \begin{matrix}
   n+\alpha+\beta+1,\, \alpha+n+1\\
    2n+\alpha+\beta+2\end{matrix}; x\Bigr]
    H(1-x)
    +\\+
    x^{-1-\beta}
 \Gamma\begin{bmatrix} 2n+\alpha+\beta+2,-\alpha-n\\
       n+1\end{bmatrix}
       \,
   \FF\Bigl[\begin{matrix} n+\alpha+\beta+1, n+1\\
              \beta+1\end{matrix};\frac 1x
	      \Bigr]
\end{multline*}
These formulae can be easily checked using the standard
technic of the Mellin--Barnes integrals, see  \cite{ASla2}
(see also tables \cite{APBM3}, 8.4.49).

We intend to evaluate
$$
\int_0^\infty K_1(x)K_2\bigl(\frac 1x\bigr)\frac {dx} x
$$
(this integral is the inner product of $\Phi_m$ and $\Phi_n$
up to a $\Gamma$-factor).

Since the Mellin transform transfer a convolution
 to a product, we must evaluate
\begin{multline*}
\frac1{2\pi i}
\int_{-i\infty}^{+i\infty}
\Gamma\begin{bmatrix} s,\,\beta+1+m-s\\
                          m+\alpha+1+s, \beta+1-s
                       \end{bmatrix}
\Gamma\begin{bmatrix}
\alpha+n+s,\, \beta+1-s\\
s,\, n+\beta+2-s
\end{bmatrix} \,ds
=
\\ =
\frac1{2\pi i}
\int_{-i\infty}^{+i\infty}
\Gamma\begin{bmatrix}
\beta+1+m-s,\, \alpha+n+s\\
m+\alpha+1+s,n+\beta+2-s
\end{bmatrix}ds
=\\
=\Gamma\begin{bmatrix}\alpha+\beta+n+m+1,\,1\\
\alpha+\beta+n+m+2,\, 1+m-n,\,1+n-m
\end{bmatrix}
\end{multline*}
(see \cite{APBM3}).
The final expression is zero if $m\ne n$
and $m-n\in\Z$.

\smallskip

{\sc Remark.} Let us explain why our calculation is not
artificial. Let us  evaluate the
 inner product
  of two Jacobi
polynomials
$$
\int_0^1
\FF\Bigl[\begin{matrix} -n,n+\alpha+\beta+1\\ \beta+1 \end{matrix};x\Bigr]
\FF\Bigl[\begin{matrix} -m,m+\alpha+\beta+1\\ \beta+1 \end{matrix};x\Bigr] x^\beta(1-x)^\alpha \,dx
$$
in the following way.
We represent the integral
as a multiplicative convolution
of
\begin{align*}
K_1(x)=
(1-x)^\alpha\FF\Bigl[\begin{matrix} -m,m+\alpha+\beta+1\\ \beta+1 \end{matrix};x\Bigr]\, H(1-x)
\\
K_2(x)=
x^{-\beta-1}\FF\Bigl[\begin{matrix} -n,n+\alpha+\beta+1\\ \beta+1 \end{matrix};\frac 1x\Bigr]\, H(x-1)
\end{align*}
and evaluate the convolution using the Mellin--Barnes
integral representations of Jacobi polynomials
(see \cite{APBM3}, 8.4).
The calculation obtained in this way admits
analytic continuation to noninteger number $n$ of
 "Jacobi polynomial" $P_n^{\alpha,\beta}$.
This gives the required formula for $K_1$, $K_2$.
\hfill $\square$

\smallskip

Orthogonal systems
 of this type appear in the following situation.
Consider the Laplace operator $\Delta$
on a pseudo-Riemannian symmetric space
$G/H$. Restricting $\Delta$ to
the subspace of $H$-invariant
functions, we obtain an ordinary second-order
differential operator. This operator is the
usual hypergeometric operator (2.3) on
some contour passing through one or more
singular points of the hypergeometric equation
(different series of symmetric spaces
give different contours). Molchanov
 obtained explicit spectral decomposition
for several boundary problems of this type
(see \cite{AMol1}, \cite{AMol2}).
These singular differential operators
have countable discrete spectra
and hence we obtain (noncomplete)
orthogonal systems consisting of
piece-wise Gauss
hypergeometric functions.
Our construction is not covered by such examples
and also does not cover them.

\smallskip

{\bf A.6. Perturbation of Wilson system.}
Koornwinder \cite{AKoo} observed that
the Wilson ${}_4F_{3}$ orthogonal system can be obtained
by application of the index hypergeometric transform
 to Jacobi polynomials.
It is possible to apply this
method to a perturbed Jacobi
system. These gives  (noncomplete)
orthogonal systems, whose elements are
pairs of ${}_4F_3$-functions.

\smallskip




{\bf A.7. An example of ${}_2F_2$-basis.} Since we discuss
possibilities of obtaining of hypergeometric bases
using standard integral transforms, let us
construct a basis that is not a deformation of a classical system.

Fix $\rho>-1$.
Consider the function $r_n$ on $\R$ given by
$$
r_n=\left\{
\begin{aligned}
(1-x)^{\rho/2} P_n^{\rho,0}(x) \qquad
  \text{for $-1<x<1$}
  \\
  0, \qquad \text{otherwise}
\end{aligned}
\right.
$$
where $P_n^{\rho,0}$ is a Jacobi polynomial.
Obviously, the functions
$r_n(x+2m)$, where $n=0,1,2,\dots$
and $m\in \Z$, form an orthogonal basis in $L^2(\R)$.

Evaluating their Fourier transforms (see \cite{APBM2}, 2.22.5.1),
we obtain the following statement

\smallskip

 {\sc Proposition A.7} {\it The functions
$$
e^{-2 m i\xi}
\,{}_2F_2\Bigl[
  \begin{matrix}\rho/2+1,\,1-\rho/2\\
                  \rho/2+n+2,-\rho/2+1-n
   \end{matrix};
    2i\xi\Bigr]
    $$
    where $n=0,1,2,\dots$ and $m\in\Z$,
    form an orthonormal basis in $L^2(\R)$.}

    \smallskip

Our expression ${}_2 F_2$ is relatively simple.
Up to a constant (depending on $n$)
it equals
$$
\sum_{k=0}^\infty \frac{(\rho/2+1+k)_n}
  {(-\rho/2+1-n+k)_n \, k!} (2i\xi)^k
  $$

  {\sc Remark.} A toy construction
   in the same type can be obtained by applying
   the Fourier transform to the standard trigonometric system
  $e^{inx}$ on $[-\pi,\pi]$. Thus the functions
  ${\sin (y\pi) e^{2i\pi my}}/  {(n+y)}$, where $m\in \Z$, $n=0,1,2,\dots$
 form an orthonormal basis in $L^2(\R)$.

\renewcommand\refname{References to Addendum}

\end{document}